\documentclass[12pt]{article}

\usepackage{amsmath,amssymb,amsfonts,amsthm,graphics,verbatim}

\begin{document}

\newcommand{\sg}{\mathbf{B_n}(\mathbb{Z}[1/p])}
\newcommand{\se}{\subseteq}
\newcommand{\rt}{\rightarrow}
\newcommand{\map}{\mathcal{A}_\mathbf{A}}
\newcommand{\ap}{\mathcal{A}}
\newcommand{\msc}{\mathfrak{S}_\mathbf{B}}
\newcommand{\sct}{\mathfrak{S}}

\title{A finitely-presented solvable group with a
 small quasi-isometry group}
\author{Kevin Wortman\thanks{Supported in part by an N.S.F. Postdoctoral
Fellowship.}}
\date{August 2, 2005}
\maketitle

\begin{abstract}

We exhibit a family of infinite, finitely-presented,
nilpotent-by-abelian groups. Each member of this family is a
solvable $S$-arithmetic group that is related to Baumslag-Solitar
groups and has a quasi-isometry group that is virtually a product
of a solvable real Lie group and a solvable $p$-adic Lie group.

In addition, we propose a candidate for a polycyclic group whose
quasi-isometry group is a solvable real Lie group, and we
introduce a candidate for a quasi-isometrically rigid solvable
group that is not finitely presented. We also record some
conjectures on the large-scale geometry of lamplighter groups.

\end{abstract}

Let $\mathbf{B_n}$ be the upper-triangular subgroup of
$\mathbf{SL_n}$, and let $\mathbf{\overline{B}_n} \leq
\mathbf{PGL_n}$ be the image of $\mathbf{B_n}$ under the natural
quotient map $\mathbf{SL_n} \rt \mathbf{PGL_n}$.

If $p$ is a prime number, then the group $\sg$ is finitely
presented for all $n \geq 2$. In particular it is finitely
generated, so we can form its quasi-isometry group---denoted
$\mathcal{QI}\big(\sg \big)$. In this paper we will prove

\bigskip \noindent \textbf{Theorem A.} \emph{If $n \geq 3$ then} $$\mathcal{QI}\big(\sg\big) \cong
\Big( \mathbf{\overline{B}_n}(\mathbb{R}) \times
\mathbf{\overline{B}_n}(\mathbb{Q}_p) \Big) \rtimes \mathbb{Z} / 2
\mathbb{Z}
$$

\bigskip

\noindent The $ \mathbb{Z} / 2 \mathbb{Z}$-action defining the
semi-direct product above is given by a $\mathbb{Q}$-automorphism
of $\mathbf{PGL_n}$ which stabilizes $\mathbf{\overline{B}_n}$.
The order $2$ automorphism acts simultaneously on each factor.

For the proof of Theorem A, the reader may advance directly to
Section 1. The proof continues in Sections 2 and 3.

\smallskip \noindent \textbf{What's new about this example.} The
group $\sg$ is solvable, and when $n \geq 3$, its quasi-isometry
group is ``small'' in two senses: it is virtually a product of
finite-dimensional Lie groups (real and $p$-adic), and it is
solvable.

\bigskip

The \emph{quasi-isometry group} of a finitely-generated group
$\Gamma$ is the group of all quasi-isometries of $\Gamma$ modulo
those that have finite distance in the sup-norm to the identity.

We don't know too much about which groups can be realized as
quasi-isometry groups, but the examples we do have suggest this
collection could be quite diverse. A fair amount of variety is
displayed even in the extremely restricted class of
finitely-generated groups that appear as lattices in semisimple
groups. These quasi-isometry groups include the $\mathbb{Q}$,
$\mathbb{R}$, or $\mathbb{Q}_p$-points of simple algebraic groups;
discrete groups that are finite extensions of lattices; and
infinite-dimensional groups (for a complete list of references to
these results see e.g. \cite{Fa} and \cite{Wo}.)

However, the same variety is not presently visible in the class of
quasi-isometry groups for infinite, finitely-generated, amenable
groups. In fact, the quasi-isometry groups of all previously
studied amenable groups contain infinite-dimensional subgroups.
Examples include infinite, finitely-generated, abelian groups (see
e.g. 3.3.B \cite{G-P}) along with many other groups from the
abundant class of abelian-by-cyclic groups.

Indeed, using the geometric models from $\cite{Fa-M Acta}$, the
quasi-isometry groups of finitely-presented, nonpolycyclic,
abelian-by-cyclic groups can be seen to contain
infinite-dimensional subgroups of automorphism groups of trees.
This generalizes the fact that $\text{Bilip}(\mathbb{Q}_m)$ is
contained in the quasi-isometry group of the Baumslag-Solitar
group $\langle a,b \mid aba^{-1}=b^m \rangle$ where
$\text{Bilip}(\mathbb{Q}_m)$ denotes the group of all bilipschitz
homeomorphisms of $\mathbb{Q}_m$. See Theorem B below.

For easy-to-describe polycyclic examples of abelian-by-cyclic
groups, take any group of the form $\Gamma = \mathbb{Z}^m \rtimes
_D \mathbb{Z}$ where $D \in \mathbf{SL_m}(\mathbb{Z})$ has a real
eigenvalue $\lambda >1$. The discrete group $\Gamma$ is a
cocompact subgroup of the Lie group $M = \mathbb{R}^m \rtimes
_\varphi \mathbb{R}$ where $\varphi (t) = D^t$. If $L_\lambda \leq
\mathbb{R}^m$ is a $k$-dimensional eigenspace of $D$ corresponding
to $\lambda$, then the cosets of $ L_\lambda \rtimes _\varphi
\mathbb{R}$ in $ M$ are the leaves of a foliation of $M$ by
totally-geodesic hyperbolic spaces. Any bilipschitz map of
$L_\lambda$ gives a quasi-isometry of each leaf and a
quasi-isometry of $M$. Since $M$ is quasi-isometric to $\Gamma$,
$\text{Bilip}(\mathbb{R}^k)$ is a subgroup of
$\mathcal{QI}(\Gamma)$. This example generalizes the well-known
fact that the quasi-isometry group of the $3$-dimensional Lie
group Sol contains $\text{Bilip}(\mathbb{R}) \times
\text{Bilip}(\mathbb{R})$.

The evidence led Farb-Mosher to a ``flexibility conjecture'': that
any infinite, finitely-generated, solvable group has an
infinite-dimensional quasi-isometry group; see $\#$3 on page 124
of \cite{F-M 2}.

Proposition A shows that this conjecture is false.

\bigskip

In addition, $\mathcal{QI}\big(\sg \big)$ is, to the best of my
knowledge, the first known example of an amenable quasi-isometry
group of an infinite, finitely-generated group.

\smallskip \noindent \textbf{Quasi-isometries of $\mathbf{S}$-arithmetic Baumslag-Solitar groups.}
The proof of Proposition A is motivated by the proof of the
following

\bigskip \noindent \textbf{Theorem B. (Farb-Mosher)}
$$\mathcal{QI}\big(\mathbf{B_2}(\mathbb{Z}[1/p])\big) \cong
\text{Bilip}(\mathbb{R}) \times \text{Bilip}(\mathbb{Q}_p)$$

\bigskip

\noindent Notice that the Baumslag-Solitar group
$\text{BS}(1,p^2)$ is an index $2$ subgroup of
$\mathbf{B_2}(\mathbb{Z}[1/p])$. (For those checking details,
$\text{Bilip}(\mathbb{Q}_p)=\text{Bilip}(\mathbb{Q}_{p^2})$.)

The above theorem is proved by studying the geometry of a space
that is essentially a ``horosphere" in a product of a hyperbolic
plane and a regular $(p+1)$-valent tree, $\mathbb{H}^2 \times T_p$
\cite{F-M 1}. Here's a brief account of the Farb-Mosher proof:

\bigskip

The group $\mathbf{SL_2}(\mathbb{R})$ acts by isometries on
$\mathbb{H}^2$, and $\mathbf{SL_2}(\mathbb{Q}_p)$ acts by
isometries on $T_p$. Since the diagonal homomorphism
$$\mathbf{B_2}(\mathbb{Z}[1/p]) \rt \mathbf{SL_2}(\mathbb{R})
\times \mathbf{SL_2}(\mathbb{Q}_p)$$ $$M \mapsto (M,M)$$ has a
discrete image, $\mathbf{B_2}(\mathbb{Z}[1/p])$ acts properly on
$\mathbb{H}^2 \times T_p$. One can show that a
$\mathbf{B_2}(\mathbb{Z}[1/p])$-orbit is a finite Hausdorff
distance from a quasi-isometrically embedded horosphere
$\mathcal{H}_p \se \mathbb{H}^2 \times T_p$. Hence, finding the
quasi-isometries of $\mathbf{B_2}(\mathbb{Z}[1/p])$ amounts to
finding the quasi-isometries of $\mathcal{H}_p$.

Farb-Mosher show that any quasi-isometry of $\mathcal{H}_p$ is the
restriction of a product of a quasi-isometry of $\mathbb{H}^2$
that preserves a foliation by horocycles, with a quasi-isometry of
$T_p$ that preserves an analogous ``foliation by horocycles". The
group of all such products can be identified with
$\text{Bilip}(\mathbb{R}) \times \text{Bilip}(\mathbb{Q}_{p})$.

Note that Theorem B can be rephrased as stating that any
quasi-isometry of the horosphere $\mathcal{H}_p \se \mathbb{H}^2
\times T_p$ is the restriction of a quasi-isometry of
$\mathbb{H}^2 \times T_p$ that stabilizes $\mathcal{H}_p$.

\bigskip

\smallskip \noindent \textbf{Relationship between Theorems A and B.}
 In the case of Theorem A,
that is when $n \geq 3$, the geometry of $\sg$ is similar to the
picture drawn by Farb-Mosher except that the rank one spaces, $
\mathbb{H}^2$ and $T_p$, are replaced by higher rank analogues. In
this case, $\sg$ is quasi-isometric to the intersection of $n-1$
distinct horospheres in the product of a higher rank symmetric
space and a higher rank Euclidean building. In comparison to the
$n=2$ case, any quasi-isometry of this intersection is the
restriction of a quasi-isometry of the entire product.

The contrast between Theorem A and Theorem B is produced by the
theorem of Kleiner-Leeb that any quasi-isometry of a product of a
higher rank symmetric space and a Euclidean building is a bounded
distance in the sup-norm from an isometry \cite{K-L}. The group of
all isometries that preserve such an intersection of horospheres
up to finite Hausdorff distance is precisely
$$\Big( \mathbf{\overline{B}_n}(\mathbb{R}) \times \mathbf{\overline{B}_n}(\mathbb{Q}_p)
\Big) \rtimes \mathbb{Z} / 2 \mathbb{Z}$$

\smallskip \noindent \textbf{Adding places to $\mathbf{S}$-arithmetic Baumslag-Solitar
groups.} Although the theme of our proof is modelled on the theme
of the Farb-Mosher proof, many of the individual techniques are
different. For example, we apply a theorem of Whyte on
quasi-isometries of coarse fibrations \cite{Wh} which tells us
that any quasi-isometry of our model space for $\sg$ preserves
what can be seen of the product decomposition between real and
$p$-adic factors.

Whyte's theorem was used in \cite{T-Wh} for the proof of the
following generalization of Theorem B.

\bigskip \noindent \textbf{Theorem C. (Taback-Whyte)} If $m$ is composite, then
$$\mathcal{QI}\big(\mathbf{B_2}(\mathbb{Z}[1/m])\big) \cong
\text{Bilip}(\mathbb{R}) \times \prod _{p |
m}\text{Bilip}(\mathbb{Q}_p)$$

\bigskip

 \bigskip

 \bigskip

\smallskip \noindent \textbf{Generalizing theorems above to conjectures below.} In what
 remains of this introduction, we state some
conjectures for groups similar to those from Theorems A, B, and C.
Two of these conjectures are new (D and F), one is well-known (E),
and the rest (G, H, and I) do not seem to be in the literature,
although they are growing in popularity amongst those who have
recognized the geometric connection between Baumslag-Solitar
groups and lamplighter groups.

All of these conjectures pertain to the large-scale geometry of
intersected horospheres in either a product of symmetric spaces,
or a product of Euclidean buildings. Whyte's theorem can not be
used in these cases to distinguish factors. This is an obstruction
to verifying the conjectures.

\smallskip \noindent \textbf{Polycyclic groups.} While Theorem
A gives an example of an infinite solvable group with a ``small"
quasi-isometry group, the question remains of whether a similar
example exists within the class of polycyclic groups.

However, Theorem A---along with its proof---provides evidence that
such an example may exist. Indeed, although $\sg$ is not
polycyclic, its geometry resembles the geometry of the polycyclic
group $\mathbf{B_n}(\mathbb{Z}[\sqrt{2}])$. The key is to note
that if $\sigma$ is the natural map derived from the nontrivial
element of the Galois group of $\mathbb{Q}(\sqrt{2}) /\mathbb{Q}$,
then the map
$$\mathbf{B_n}(\mathbb{Z}[\sqrt{2}]) \rt
\mathbf{SL_n}(\mathbb{R})\times \mathbf{SL_n}(\mathbb{R})$$
defined by $$M \mapsto (M,\sigma(M))$$ is an injective
homomorphism with discrete image.

The polycyclic group above is quasi-isometric to an intersection
of $n-1$ horospheres in the symmetric space for
$\mathbf{SL_n}(\mathbb{R}) \times \mathbf{SL_n}(\mathbb{R})$.
Asking that every quasi-isometry of the intersection be the
restriction of a quasi-isometry of the symmetric space is the
content of

\bigskip \noindent \textbf{Conjecture D.} \emph{If $n \geq 3$, then}
$$\mathcal{QI}\big( \mathbf{B_n}(\mathbb{Z}[\sqrt{2}]) \big) \cong \Big( \mathbf{\overline{B}_n}(\mathbb{R}) \times
\mathbf{\overline{B}_n}(\mathbb{R}) \Big)
 \rtimes \Big( \mathbb{Z} / 2 \mathbb{Z} \times \mathbb{Z} / 2 \mathbb{Z} \Big) $$

\bigskip

Note that Theorem A is to Theorem B as Conjecture D is to the
well-known

\bigskip \noindent \textbf{Conjecture E.}
$$\mathcal{QI}\big(\mathbf{B_2}(\mathbb{Z}[\sqrt{2}])\big) \cong
\Big( \text{Bilip}(\mathbb{R}) \times
\text{Bilip}(\mathbb{R})\Big) \rtimes \mathbb{Z} / 2\mathbb{Z}$$

The above is better-known as a conjecture about the
$3$-dimensional real Lie group Sol---the geometry associated with
a $3$-manifold that fibers over a circle with torus fibers and
Anosov monodromy. Indeed, $\mathbf{B_2}(\mathbb{Z}(\sqrt{2}))$ is
a discrete cocompact subgroup of Sol. (See \cite{F-M 2} for
remarks concerning Conjecture E.)

There is almost nothing known about quasi-isometry types of
polycyclic (and not virtually nilpotent) finitely-generated
groups, and Conjecture E stands as one of the more important open
questions in geometric group theory in part because it is seen as
a first step towards understanding the large-scale geometry of
polycyclic groups.

\smallskip \noindent \textbf{Geometric comparisons between rank one amenable groups.}
 Note that Sol is quasi-isometric to a
horosphere in $\mathbb{H}^2 \times \mathbb{H}^2$, and Conjecture E
asserts that any quasi-isometry of Sol extends to a product of
quasi-isometries of $\mathbb{H}^2$. Each of the quasi-isometries
of $\mathbb{H}^2$ would preserve a foliation by horocycles, and
the supremum of all Hausdorff distances between horocycles and
their images would be finite.

This is in analogy with the known large-scale geometry of solvable
Baumslag-Solitar groups described above, and in analogy with a
conjectural picture of lamplighter groups that is discussed below.
Each of these three examples of groups are rank one amenable
groups.

\smallskip \noindent \textbf{Function-field-arithmetic groups.}
Theorems A, B, and C are statements about $S$-arithmetic groups,
while Conjectures D and E concern arithmetic groups. The geometry
of these groups extends easily to the geometry of
function-field-arithmetic groups.

Specifically, let $\mathbb{F}_q((t))$ be a field of Laurent series
with coefficients in a finite field $\mathbb{F}_q$, and let
$\mathbb{F}_q[t,t^{-1}]$ be the ring of Laurent polynomials.
Expanding on Theorem A and Conjecture D we have the following

\bigskip \noindent \textbf{Conjecture F.} \emph{If $n \geq 3$, then}
$$\mathcal{QI}\big( \mathbf{B_n}(\mathbb{F}_q[t,t^{-1}]) \big)
 \cong \Big( \mathbf{\overline{B}_n}\big(\mathbb{F}_q((t))\big) \times
\mathbf{\overline{B}_n}\big(\mathbb{F}_q((t))\big) \Big)
 \rtimes  \Big( \mathbb{Z} / 2 \mathbb{Z} \times \mathbb{Z} / 2 \mathbb{Z} \Big) $$

\bigskip

\noindent Bux showed that the group
$\mathbf{B_n}(\mathbb{F}_q[t,t^{-1}])$ is finitely generated but
not finitely presented \cite{Bu}. This group acts on a product of
higher rank buildings.

Turning our attention to the case $n=2$, the group
$\mathbf{B_2}(\mathbb{F}_q[t,t^{-1}])$ acts on a horosphere in a
product of $(q+1)$-valent trees. Theorem B and Conjecture E lead
to

\bigskip \noindent \textbf{Conjecture G.}
$$\mathcal{QI}\big(\mathbf{B_2}(\mathbb{F}_q[t,t^{-1}])\big) \cong
\Big( \text{Bilip}\big(\mathbb{F}_q((t))\big) \times
\text{Bilip}\big(\mathbb{F}_q((t))\big)\Big) \rtimes \mathbb{Z} /
2\mathbb{Z}$$

\bigskip

\noindent The groups $\mathbf{B_2}(\mathbb{F}_q[t,t^{-1}])$ are
also finitely generated and not finitely presented.

\smallskip \noindent \textbf{Lamplighter groups.} Conjecture G is closely
related to the study of quasi-isometry types of lamplighter
groups.

Recall that if $G$ is a finite group, then the wreath product
$$G \wr \mathbb{Z} = \Big(\bigoplus _{\mathbb{Z}}G \Big) \rtimes \mathbb{Z}$$
is called a \emph{lamplighter group}.

Note that the group $\mathbf{B_2}(\mathbb{F}_2[t,t^{-1}])$ is
isomorphic to the wreath product $\mathbb{F}_2^2 \wr \mathbb{Z}$.
This example reveals how the natural geometry of the
function-field-arithmetic groups from Conjecture G can be formally
generalized---in the same way that Baumslag-Solitar groups can be
viewed as a formal generalization of the $S$-arithmetic groups
$\mathbf{B_2}(\mathbb{Z}[1/p])$---to act on a horosphere in a
product of trees. These horospheres are examples of
\emph{Diestel-Leader graphs}.

In Section 4, we will briefly trace through a description of the
geometric action of lamplighters on a horosphere in a product of
trees. The analogy with solvable Baumslag-Solitar groups and with
Sol leads to the following two conjectures.

\bigskip \noindent \textbf{Conjecture H.} \emph{For any nontrivial finite group $G$,}
$$\mathcal{QI}\big(G \wr \mathbb{Z} \big) \cong
\Big( \text{Bilip}\big(G((t))\big) \times
\text{Bilip}\big(G((t))\big)\Big) \rtimes \mathbb{Z} /
2\mathbb{Z}$$

\bigskip

\bigskip \noindent \textbf{Conjecture I.} \emph{Let $G$ and $H$ be finite groups.
Then $G \wr \mathbb{Z}$ is quasi-isometric to $H \wr \mathbb{Z}$
if and only if $|G|^k=|H|^j$ for some $k,j \in \mathbb{N}$.}

\bigskip

The question of when two lamplighter groups are quasi-isometric
was asked by de$\,$la$\,$Harpe; see IV.B.44 \cite{Ha}. The ``if"
implication in the above conjecture is well-known (and easily
seen) to be true, but there does not exist a single known case of
the ``only if" implication. For example, it is unknown even if
$(\mathbb{Z}/2\mathbb{Z} ) \wr \mathbb{Z}$ is quasi-isometric to
$(\mathbb{Z}/3\mathbb{Z} ) \wr \mathbb{Z}$.

Note the similarity between the above conjecture and

\bigskip \noindent \textbf{Theorem J. (Farb-Mosher)} \emph{If $\text{\emph{BS}}(1,m)$
 denotes the Baumslag-Solitar group $\langle a,b \mid
aba^{-1}=b^{m} \rangle$, then $\text{\emph{BS}}(1,m)$ is
quasi-isometric to $\text{\emph{BS}}(1,r)$ if and only if
$m^k=r^j$ for some $k,j \in \mathbb{N}$.}

\bigskip

\noindent See Theorem 7.1 of \cite{F-M 1}.

\bigskip

\smallskip \noindent \textbf{Lattices in solvable groups.} With the exception of the
general lamplighter groups $G \wr \mathbb{Z}$ and the general
Baumslag-Solitar groups $\text{BS}(1,m)$, all of the
finitely-generated groups mentioned above are discrete cocompact
subgroups of solvable Lie groups over locally-compact, nondiscrete
fields. This follows from reduction theory.

\bigskip \noindent \textbf{Realizing quasi-isometry groups.} As
mentioned above, a motivation for writing this paper was to add to
the list of examples of groups that are known to be realized as
quasi-isometry groups of finitely-generated groups. It would be
nice to know more examples.

\smallskip \noindent \textbf{Notation.} In the remainder of this paper we
let
$$G_p = \Big( \mathbf{\overline{B}_n}(\mathbb{R}) \times \mathbf{\overline{B}_n}(\mathbb{Q}_p)
\Big) \rtimes \mathbb{Z} / 2 \mathbb{Z}$$

\bigskip \noindent \textbf{Outline of paper.} Our proof of Proposition A
 begins in Section 1 where we describe
 a geometric model for $\sg$. We use the
model in Section 2 to show that $G_p$ is included in the
quasi-isometry group of $\sg$, and we use the model in Section 3
to show that every quasi-isometry of $\sg$ is given by an element
of $G_p$.

Section 1 does not use the assumption that $n \geq 3$, and Section
2 may also be read assuming that $n=2$ as long as one replaces
$G_p$ with its obvious index 2 subgroup. Even Section 3 only
invokes the $n\geq 3$ assumption in efficient and isolated
applications of Kleiner-Leeb's theorem. We point all of this out
because the reader may find it helpful to fall back on the case
$n=2$ at times. Indeed, it's important to have a solid
understanding of the model space while reading the first three
sections, and it is possible to draw and easiest to visualize the
model space when $n=2$. For descriptions of the model space when
$n=2$, see e.g. 7.4 of \cite{word processing} or \cite{F-M 1}. The
definitions of the model space in these references is different
than the one we will give, but they are easily seen to be
equivalent.

Section 4 contains a sketch of the geometry of lamplighters that
hints at Conjectures H and I.

\bigskip

\smallskip \noindent \textbf{Acknowledgements.} In the nearly two years it took me
 to begin writing this paper after having first sketched the proof of Theorem A, I was happy
 to discuss the contents of this paper with many of my colleagues.

 Benson
 Farb convinced me to write this paper. I thank him for that,
  for his continued support, for explaining to me the context of the results in this paper,
   and for comments made on an earlier draft.

  Thanks to Tullia Dymarz and David Fisher for comments on an
earlier draft, and to Lee Mosher for two very helpful
conversations we had a couple of years ago regarding the
large-scale geometry of solvable groups. I am also appreciative of
conversations that I had with Jennifer Taback and Sean Cleary on
the contents of Section 4 and on Diestel-Leader graphs. I thank
Kevin Whyte for explaining to me that his results could be used in
the proof below.

Several of my colleagues helped to shape the material in this
paper in other ways: Angela Barnhill, Tara Brendle, Nathan
Broaddus, Indira Chatterji, Alex Eskin, Paul Jung, Irine Peng,
Alireza Salehi-Golsefidy, Edward Swartz, Karen Vogtmann, and
Wolfgang Woess. I value their contribution.

I thank Dan Margalit for his encouragement.


\pagebreak

 \noindent \Large \textbf{1. The model space for $\sg$} \normalsize
\bigskip

Throughout, we fix $n\geq 3$. Our goal in this section is to
define a metric space $X$ that is quasi-isometric to $\sg$. Then
$\mathcal{QI}(X) \cong \mathcal{QI}\big(\sg \big) $. In Sections 2
and 3, we will determine $\mathcal{QI}( X)$.

As alluded to in the introduction, $X$ will essentially be an
intersection of $n-1$ horospheres in a product of a symmetric
space and a Euclidean building, but we will define it differently
to simplify some points of our proof.

Although we are assuming $n \geq 3$ in Theorem A, the entire
contents of this section hold with the weaker assumption that $n
\geq 2$.

\bigskip \noindent \textbf{Quasi-isometries.}
 For constants $K \geq 1$ and $C \geq 0$, a $(K,C)$
\emph{quasi-isometric embedding} of a metric space $Y$ into a
metric space $Z$ is a function $\phi : Y \rightarrow Z$ such that
for any $y_{1},y_{2} \in Z$:
$$\frac{1}{K} d\big(y_{1},y_{2}\big) -C \leq d\big(\phi(y_{1}),\phi(y_{2})\big)
\leq   K d\big(y_{1},y_{2}\big) +C$$

We call $\phi$ a $(K,C)$ \emph{quasi-isometry} if $\phi$ is a
$(K,C)$ quasi-isometric embedding and there is a number $D\geq 0$
such that every point in $Z$ is within distance $D$ of some point
in the image of $Y$.

\bigskip \noindent \textbf{Quasi-isometry groups.} For a metric space $Y$, we define the relation $\sim$  on the set
of functions  $Y\rt Y$ by $\phi \sim \psi$ if
$$\sup _{y\in Y} d\big(\phi(y),\psi(y)\big)<\infty$$

 We denote
the quotient space of all quasi-isometries of $Y$ modulo $\sim$ by
$\mathcal{QI}( Y )$. We call $\mathcal{QI}( Y )$ the
\emph{quasi-isometry group} of $Y$ as it has a natural group
structure arising from function composition.

If $\Gamma$ is a finitely-generated group, we use its
left-invariant word metric to form the quasi-isometry group
$\mathcal{QI}( \Gamma )$.

\medskip \noindent \large \textbf{1.1 Description of the model space} \normalsize \medskip

The finitely-presented group $\sg$ is a discrete subgroup of
$\mathbf{SL_n}(\mathbb{R}) \times \mathbf{SL_n}(\mathbb{Q}_p)$
and, therefore, acts on a product of nonpositively curved spaces.

\smallskip \noindent \textbf{Symmetric space and Euclidean building.}
 We let $X_\infty$ be the symmetric
space for $\mathbf{SL_n}(\mathbb{R})$, and $X_p$ be the Euclidean
building for $\mathbf{SL_n}(\mathbb{Q}_p)$. All of the facts we
will use about $X_\infty$ and $X_p$ are standard and can be found,
for example, in \cite{B-G-S} or \cite{Br}.

\smallskip \noindent \textbf{Topological description of $\mathbf{X}$.} We let $\mathbf{U_n}$ be the subgroup of $\mathbf{B_n}$ consisting
of matrices whose diagonal entries are all equal to $1$. We define
$X$ to be the topological space $X_p \times
\mathbf{U_n}(\mathbb{R})$ and we let $\pi : X \rt X_p$ be the
projection map.

\smallskip \noindent \textbf{A symmetric space in $\mathbf{X}$.}
Denote the diagonal subgroup of $\mathbf{B_n}$ by $\mathbf{A_n}$.
We let $\map \se X_p$ be the apartment corresponding to
$\mathbf{A_n}(\mathbb{Q}_p)$ and $F_\mathbf{A} \se X_\infty$ be
the flat corresponding to $\mathbf{A_n}(\mathbb{R})$. Recall that
with the metric restricted from $X_p$ and $X_\infty$ respectively,
$\map$ and $F_\mathbf{A}$ are each isometric to
$(n-1)$-dimensional Euclidean space.

The group of diagonal matrices $ \mathbf{A_n}(\mathbb{Z}[1/p])$
acts on the Euclidean spaces $\map$ and $F_\mathbf{A}$ as a
discrete group of translations. These two isometric actions are
related:

\bigskip \noindent \textbf{Lemma 1.1(a).} \emph{After possibly scaling
 the metric on $X_p$, there exists an $
\mathbf{A_n}(\mathbb{Z}[1/p])$-equivariant isometry}
 $$f : \map \rt F_\mathbf{A}$$

\bigskip

\noindent \textbf{Proof.} Choose a sector $\msc \se \map$ that is
fixed up to finite Hausdorff distance by the action of
$\mathbf{B_n}(\mathbb{Q}_p)$ on $X_p$, and let
$\mathfrak{C}_\mathbf{B} \se F_\mathbf{A}$ be a Weyl chamber that
is fixed up to finite Hausdorff distance by the action of
$\mathbf{B_n}(\mathbb{R})$ on $X_\infty$.

We denote the collection of roots of $\mathbf{SL_n}$ with respect
to $\mathbf{A_n}$ that are positive in $\mathbf{B_n}$ as $\alpha
_1 , \alpha _2 ,..., \alpha _{n-1}$. Each $\alpha _i$ corresponds
to a wall of $\msc$ (respectively $\mathfrak{C}_\mathbf{B}$) which
we denote by $L_{p,i}$ (respectively $L_{\infty, i}$) for all $i$,
so there is a unique isometry that maps $\msc$ onto
$\mathfrak{C}_\mathbf{B}$ and $L_{p,i}$ onto $L_{\infty, i}$ for
each $i$. We name this isometry $f':\map \rt F_\mathbf{A}$.

We denote the inverse transpose automorphism of $\mathbf{SL_n}$ by
$\varphi ^{\iota t}$, and let $\varphi ^{\iota t} _\infty$ be the
isometry of $X_\infty$ induced by $\varphi ^{\iota t}$. Notice
that $\varphi ^{\iota t} _\infty (F_\mathbf{A})=F_\mathbf{A}$, so
the composition $\varphi ^{\iota t} _\infty \circ f'$ defines an
isometry $f:\map \rt F_\mathbf{A}$.

For $1 \leq i \leq n-1$, we define the diagonal matrix $a_i \in
\mathbf{A_n}(\mathbb{Z}[1/p])$ as

\[
(a_i)_{\ell,k}=
\begin{cases}
p, & \text{if $\ell = k = i$;}\\
1/p, & \text{if $\ell = k = i+1$;}\\
\delta_{\ell,k}, & \text{otherwise.}
\end{cases}
\]

Let $\beta$ be any root of $\mathbf{SL_n}$ with respect to
$\mathbf{A_n}$, and let $\mathbf{U}_\beta \leq \mathbf{SL_n}$ be
the corresponding root group. If $a_i$ acts as an expanding
automorphism on $\mathbf{U}_\beta (\mathbb{Q}_p)$, then it acts as
a contracting automorphism on $\mathbf{U}_\beta (\mathbb{R})$
since the $p$-adic norm of $p$ is less than $1$. Thus, if $a_i$
translates $\map$ towards the chamber at infinity corresponding to
$\mathbf{U}_{\beta}$, then $a_i$ translates $F_\mathbf{A}$ towards
the chamber at infinity corresponding to
$\mathbf{U}_{-\beta}=\varphi ^{\iota t}(\mathbf{U}_\beta )$.
Therefore, after rescaling $X_p$ so that the translation lengths
of each of the $a_i$ are the same on $\map$ as they are on
$F_\mathbf{A}$, we have that $f$ is $\langle a_1 , a_2, ...,
a_{n-1} \rangle $-equivariant.

Our lemma follows since the
$\mathbf{A_n}(\mathbb{Z}[1/p])$-actions factor through the group
$\langle a_1 , a_2, ..., a_{n-1} \rangle $.

\hfill $\blacksquare$

We write the connected component of the identity in
$\mathbf{A_n}(\mathbb{R})$ as $\mathbf{A_n}(\mathbb{R})^\circ$,
and recall that it acts transitively on $F_\mathbf{A}$ without
fixed points. Thus, there is a natural bijection between
$F_\mathbf{A}$ and $\mathbf{A_n}(\mathbb{R})^\circ$.

With the Weyl chamber $\mathfrak{C}_\mathbf{B} \se F_\mathbf{A}$
as in the proof of Lemma 1.1(a), we have:
\begin{quote} (\emph{i}) every point in $X_\infty$ is an element
of a maximal flat whose $R$-neighborhoods contain
$\mathfrak{C}_\mathbf{B}$ when $R \gg 0$;\\
(\emph{ii}) the group $\mathbf{U_n}(\mathbb{R})$ acts transitively
on the set of all maximal flats as in (\emph{i}); \\
(\emph{iii}) $F_\mathbf{A}$ is a flat as in (\emph{i}); and \\
(\emph{iv}) $\mathbf{A_n}(\mathbb{R})^\circ$ acts transitively on
$F_\mathbf{A}$. \end{quote} Therefore, the group
$\mathbf{A_n}(\mathbb{R})^\circ\mathbf{U_n}(\mathbb{R})$ acts
transitively on $X_\infty$. It is also easy to see that this
action is without fixed points. Thus, there is an isometry (up to
scale) between $X_\infty$ and the Lie group
$\mathbf{B_n}(\mathbb{R})^\circ =
\mathbf{A_n}(\mathbb{R})^\circ\mathbf{U_n}(\mathbb{R})$ endowed
with a left-invariant metric. This allows us to identify
$X_\infty$ with $F_\mathbf{A} \times \mathbf{U_n}(\mathbb{R})$.

Note that $ \pi ^{-1} (\map) = \map \times
\mathbf{U_n}(\mathbb{R})$, so we define the diffeomorphism
$$\widehat{f}: \pi ^{-1} (\map) \rt X_\infty$$ by
$$\widehat{f}(a,u)=(f(a),u)$$ for all $a\in \map$ and all $u \in
\mathbf{U_n}(\mathbb{R})$.

If $\omega$ is the metric on
$X_\infty$, we endow $\pi ^{-1} (\map)$ with the pull-back metric
$\widehat{f}^*(\omega)$.

\smallskip \noindent \textbf{A family of symmetric spaces in $\mathbf{X}$.}
Recall from the proof of Lemma 1.1(a) that we chose a sector $\msc
\se \map$ that is fixed up to finite Hausdorff distance by the
action of $\mathbf{B_n}(\mathbb{Q}_p)$ on $X_p$. We call any
apartment in $X_p$ that contains a subsector of $\msc$ a
\emph{symmetric apartment}. If $\ap \se X_p$ is a symmetric
apartment, then we call $\pi ^{-1}(\ap )$ a \emph{symmetric
pre-image}.

Denote the building retraction corresponding to the pair $(\map ,
\msc)$ by $$\varrho _{\map , \msc} : X_p \rt \map$$ We extend
$\varrho _{\map , \msc}$ to a retraction of the entire space $X$
by defining
$$\widehat{\varrho} _{\map , \msc} : X \rt \pi ^{-1} ( \map )$$
as $$\widehat{\varrho} _{\map , \msc} \big(x,u\big) = \big(
\varrho _{\map , \msc}(x),u\big)$$ for all $x \in X_p$ and $u \in
\mathbf{U_n}(\mathbb{R})$.

Since $\varrho _{\map , \msc}$ maps symmetric apartments
isometrically onto $\map$, the map $\widehat{\varrho} _{\map ,
\msc}$ defines a diffeomorphism of any symmetric pre-image onto
$\pi ^{-1} ( \map )$. We endow each symmetric pre-image with the
pull-back metric $(f \circ \widehat{\varrho} _{\map ,
\msc})^*(\omega)$.

With the metric $(f \circ \widehat{\varrho} _{\map ,
\msc})^*(\omega)$, each symmetric pre-image is isometric to the
symmetric space $X_\infty$. Also note that the metrics on any pair
of symmetric pre-images agree on their intersection.

\smallskip \noindent \textbf{A metric on $\mathbf{X}$.}
Recall that any point in $X_p$ is contained in a symmetric
apartment. Therefore, every point in $X$ is contained in a
symmetric pre-image. This allows us to endow $X$ with the path
metric.

The two following lemmas illustrate how the metric space $X$
reflects the geometry of the Euclidean building and the symmetric
pre-images that defined it.

\bigskip \noindent \textbf{Lemma 1.1(b).} \emph{The restricted metric from $X$
 and the metric $(f \circ \widehat{\varrho} _{\map , \msc})^*(\omega)$ are equal on any symmetric pre-image.}

\bigskip

\noindent \textbf{Proof.} If $\gamma$ is a path in $X$, then
$\widehat{\varrho} _{\map , \msc}(\gamma)$ is a path in $X$ that
is no longer than $\gamma$. This shows that $\pi ^{-1} ( \map )$
is isometric to $X_\infty$, and a similar argument applies to any
other symmetric pre-image.

\hfill $\blacksquare$

Let $\sigma : X_p \rt X$ be the section of $\pi :X \rt X_p$
defined by $\sigma (x) = (x,1)$.

\bigskip \noindent \textbf{Lemma 1.1(c).} \emph{The section
$\sigma : X_p \rt X$ is an isometric inclusion.}

\bigskip

Before proving Lemma 1.1(c), we will examine a metric property of
$\pi$.

\smallskip \noindent \textbf{The projection $\pi$ is distance nonincreasing.}
For any symmetric apartment $\mathcal{A}$, we let $\pi
_\mathcal{A}$ be the restriction of $\pi $ to $\pi ^{-1}
(\mathcal{A} )$.

Notice that if $\pi _\mathcal{A}$ is distance nonincreasing for
all $\mathcal{A}$, then $\pi$ is distance nonincreasing since any
path in $X$ can be written as a finite union of paths contained in
a symmetric pre-image.

Furthermore, if $\pi _{\map}$ is distance nonincreasing, then
every $\pi _\mathcal{A}$ is distance nonincreasing since
$$\pi_{\map} \circ \widehat{\varrho} _{\map , \msc} =
\widehat{\varrho} _{\map , \msc} \circ \pi _\mathcal{A} $$ So to
check that $\pi$ is distance nonincreasing, we only need to check
that $\pi _{\map}$ is. In other words, we want to show that the
map $$q: \mathbf{B_n}(\mathbb{R})^\circ \rt
\mathbf{A_n}(\mathbb{R})^\circ$$ $$q(au)=a$$ is distance
nonincreasing.

We proceed by letting $v$ be a tangent vector to
$\mathbf{B_n}(\mathbb{R})^\circ$ at a point $au \in
\mathbf{B_n}(\mathbb{R})^\circ$. We let $u'=au^{-1}a^{-1}$, we
denote the derivative of a map $g$ by $g_*$, and we denote the
Riemannian norm of a vector by $\| \cdot \|$. Then, $$\| v \|
_{au} =\| u'_* v \| _a \geq \| q_* u'_* v \| _a $$

It is easy to check that $q_* u'_* v  =q_* v $ since $u' \in
\mathbf{U_n}(\mathbb{R})$. Therefore, $$\| v \|_{au} \geq \| q_* v
\| _a $$ It follows that $q$, and hence $\pi$, is distance
nonincreasing.

\bigskip

We are ready to return to the

\bigskip \noindent \textbf{Proof of 1.1(c).} Clearly $\sigma$ is
distance nonincreasing, so we only need to show that if $x,y \in
X_p$, then
$$d(x,y)\leq d(\sigma (x), \sigma (y))$$
But if $\gamma \se X$ is a path between $\sigma (x)$ and $\sigma
(y)$, then $\pi (\gamma )$ is a path between $x$ and $y$ and thus
$$d(x,y) \leq \text{length}(\pi (\gamma)) \leq \text{length}( \gamma)$$

\hfill $\blacksquare$

\medskip \noindent \large \textbf{1.2 $\mathbf{X}$ is a model space} \normalsize \medskip

We will use the Milnor-\v{S}varc lemma to show that $\sg$ and $X$
are quasi-isometric. The next three lemmas are meant to verify the
hypothesis of that lemma.

\bigskip \noindent \textbf{Lemma 1.2(a).} \emph{There is an isometric $\sg$-action on $X$.}

\bigskip

\noindent \textbf{Proof.} For any $x\in X_p$, we let $a_x \in
\mathbf{A_n}(\mathbb{R})^\circ$ be the group element that is
identified with $f \circ \varrho _{\map , \msc}(x) \in
F_\mathbf{A}$ via the action of $\mathbf{A_n}(\mathbb{R})^\circ$
on $F_\mathbf{A}$.

For $b \in \sg$, we let $a_b \in \mathbf{A_n}(\mathbb{Z}[1/p])$
and $u_b \in \mathbf{U_n}(\mathbb{Z}[1/p])$ be such that
$b=a_bu_b$.

Since $\mathbf{B_n}(\mathbb{R})$ acts by isometries on $X_\infty$,
the group $\sg$ acts by isometries on $\pi ^{-1}(\map)$. With our
identifications made earlier, this action is given by
$$b*(x,u)=(a_bx, a_x^{-1}u_ba_xu)$$
for all $x \in \map$ and all $u \in \mathbf{U_n}(\mathbb{R})$.

The above action is not proper, and it is of no use to us aside
from motivating our desired action of $\sg$ on $X$ that is defined
by
$$b(x,u)=(bx, a_x^{-1}u_ba_xu)$$
for all $(x,u) \in X_p \times \mathbf{U_n}(\mathbb{R})$.

Any $b \in \sg$ restricts to an isometry $\pi ^{-1}(\mathcal{A})
\rt \pi ^{-1}(b\mathcal{A})$ for each symmetric pre-image
$\mathcal{A}$. It follows that $\sg$ acts by isometries on $X$ as
any path in $X$ can be written as a finite union of paths that are
each contained in a symmetric pre-image.

\hfill $\blacksquare$

\bigskip \noindent \textbf{Lemma 1.2(b).} \emph{The $\sg$-action on $X$ is proper.}

\bigskip

\noindent \textbf{Proof.} From the description of the metric space
$X$, it is clear that if $\{b_n \}_{n\in \mathbb{N}}$ is a
sequence of elements in $\sg$ that converges to $1$ in
$\text{Isom}(X)$, then $b_n \rightarrow 1$ in $\text{Isom}(X_p)$
and in $\text{Isom}(X_\infty)$. Thus, $b_n \in \{\pm \text{Id}\}$
for all $n \gg 0$.

\hfill $\blacksquare$

\bigskip \noindent \textbf{Lemma 1.2(c).} \emph{The $\sg$-action on $X$ is cocompact.}

\bigskip

\noindent \textbf{Proof.} Let the vertex $v \in \map$ be the
unique point in $X_p$ that is fixed by the action of
$\mathbf{SL_n}(\mathbb{Z}_p)$, and let $\mathfrak{S}_v \se \map$
be the unique sector that contains $v$ as its cone-point and such
that $\mathfrak{S}_v \cap \msc$ contains a sector. Note that
$\mathbf{U_n}(\mathbb{Z})$ fixes every point in $\mathfrak{S}_v$.

Since $\mathbf{A_n}(\mathbb{Z}[1/p])$ acts cocompactly on $\map$,
there is a compact subset $C_p \se \mathfrak{S}_v$ such that
$\mathbf{A_n}(\mathbb{Z}[1/p]) C_p =\map$.

Because $\mathbf{U_n}(\mathbb{Z}[1/p])$ is a dense subgroup of
$\mathbf{U_n}(\mathbb{Q}_p)$, and since the latter group acts
transitively on the set of symmetric apartments, we have that
\begin{align*} \mathbf{B_n}(\mathbb{Z}[1/p]) C_p & = \mathbf{U_n}(\mathbb{Z}[1/p]) \mathbf{A_n}(\mathbb{Z}[1/p]) C_p\\
 & = \mathbf{U_n}(\mathbb{Z}[1/p]) \map \\
 & = X_p \end{align*}

It is well-known that the $\mathbf{U_n}(\mathbb{Z})$-action on
$\mathbf{U_n}(\mathbb{R})$ has a compact fundamental domain which
we name $C_\infty$. We use $C_\infty$ to define a compact subset
of $X$:
$$C=\bigcup _{x\in C_p}\Big( \{x\} \times a^{-1}_x C_\infty a_x \Big)$$

Now \begin{align*} \sg C & = \sg
\mathbf{U_n}(\mathbb{Z})C  \\
& = \sg \bigcup _{x \in C_p}\Big( \{ x \}
 \times a^{-1} _x \mathbf{U_n}(\mathbb{Z}) a_x a^{-1}_x C_\infty a_x \Big) \\
& = \sg \bigcup _{x\in C_p} \Big( \{ x \} \times \mathbf{U_n}(\mathbb{R}) \Big) \\
& = X
\end{align*}

\hfill $\blacksquare$

Combining the three lemmas above yields

\bigskip \noindent \textbf{Lemma 1.2(d).} \emph{There is an isomorphism of groups}
 $$\mathcal{QI}\big(\sg\big)\cong \mathcal{QI}\big(X\big)$$

\bigskip

\noindent \textbf{Proof.} Apply the Milnor-\v{S}varc lemma.

\hfill $\blacksquare$


\medskip \noindent \large \textbf{1.3 Aside on finiteness properties.} \normalsize \medskip

Since $X$ is contractible, $\sg$ is of type $F_\infty$.

\pagebreak

 \noindent \Large \textbf{2. Some quasi-isometries of
the model space} \normalsize \bigskip

Recall that we let $G_p$ be the group $$\Big(
\mathbf{\overline{B}_n}(\mathbb{R}) \times
\mathbf{\overline{B}_n}(\mathbb{Q}_p) \Big) \rtimes \mathbb{Z} / 2
\mathbb{Z}$$ The goal of this section is to show that $G_p \leq
\mathcal{QI}(X)$.

In the case when $n=2$, one can replace all occurrences of the map
$\varphi ^*$ below with the identity to prove that the obvious
index 2 subgroup of $G_p$ is contained in $\mathcal{QI}(X)$.

\smallskip \noindent \textbf{A non-inner automorphism of $\mathbf{PGL_n}$.}
Because $n \geq 3$, there is an order $2$ automorphism of
$\mathbf{PGL_n}$ that is defined over $\mathbb{Q}$, stabilizes
$\mathbf{\overline{B}_n}$, and is not type preserving on the
spherical building for $\mathbf{PGL_n}(\mathbb{Q})$. Name this
automorphism $\varphi ^*$.

\smallskip \noindent \textbf{Isometries introduced from the $\mathbf{p}$-adic base.}
The group of isometries of $X_p$ is identified with
$\mathbf{PGL_n}(\mathbb{Q}_p) \rtimes \langle \varphi^{*}
\rangle$. The subgroup of isometries that map a subsector of
$\msc$ to another subsector of $\msc$ is identified with
$\mathbf{\overline{B}_n}(\mathbb{Q}_p) \rtimes \langle \varphi^{*}
\rangle$.

For any $\alpha \in \mathbf{\overline{B}_n}(\mathbb{Q}_p) $, we
define $$H_p(\alpha ) : X \rt X$$ by $H_p(\alpha )(x,u)=(\alpha
(x), u)$. We claim that $H_p(\alpha )$ is an isometry.

Indeed, if $\ell \in \mathbf{\overline{B}_n}(\mathbb{Q}_p)$ fixes
a subsector of $\msc$ pointwise, then $H_p(\ell)$ simply permutes
the symmetric pre-images via diffeomorphisms that preserve the
metrics on the symmetric pre-images as can easily be seen from the
definition of the metrics.

Alternatively, if $\tau \in \mathbf{\overline{B}_n}(\mathbb{Q}_p)$
restricts to a translation of $\map$, then it follows from our
identification of $\pi ^{-1} (\map )$ with
$\mathbf{A_n}(\mathbb{R})^\circ \mathbf{U_n}(\mathbb{R})$ that
$H_p(\tau )|_{\pi ^{-1}(\map)}$ is an isometry. Therefore,
$H_p(\tau )$ is also an isometry when restricted to any other
symmetric pre-image $\pi ^{-1}(\mathcal{A})$  since
$$H_p(\tau)|_{\pi ^{-1}(\mathcal{A})} = H_p(\tau  \ell ^{-1}  \tau ^{-1} )
\circ H_p (\tau )|_{\pi ^{-1}(\map)} \circ H_p (\ell ) |_{ \pi
^{-1}(\mathcal{A})}$$ where $\ell \in
\mathbf{\overline{B}_n}(\mathbb{Q}_p)$ is as in the above
paragraph with $\ell \mathcal{A} = \map$. (We used here that $\tau
\ell ^{-1} \tau ^{-1}$ fixes a subsector of $\msc $ pointwise.) It
follows that $H_p(\tau )$ is an isometry.

Our claim is substantiated since any $\alpha \in
\mathbf{\overline{B}_n}(\mathbb{Q}_p)$ is a composition of some
$\ell$ and $\tau$ as above. Thus, we have defined a homomorphism
 $$H_p :
\mathbf{\overline{B}_n}(\mathbb{Q}_p) \rt \mathcal{QI}(X)$$

\smallskip \noindent \textbf{Bilipschitz maps introduced from the symmetric pre-images.}
 The isometry group of $\pi ^{-1}(\map)$ (and of any symmetric pre-image) is identified with
$\mathbf{PGL_n}(\mathbb{R}) \rtimes \langle \varphi ^* \rangle$.

Let $\msc ^{\text{op}} \se \map$ be a sector opposite to $\msc$.
Then $\sigma ( \msc ^{\text{op}} )$ is a Weyl chamber in $\pi
^{-1}(\map)$, and with our definition of $f$ from 1.1(a), $\sigma
(\msc ^{\text{op}})$ equals $\widehat{f}^{-1}(\mathfrak{C
_\mathbf{B} })$ up to a translation of $\sigma( \map )$. Hence,
the subgroup of those isometries in $\mathbf{PGL_n}(\mathbb{R})
\rtimes \langle \varphi ^* \rangle$ that stabilizes $\sigma ( \msc
^{\text{op}} )$ up to a finite Hausdorff distance is precisely
$\mathbf{\overline{B}_n}(\mathbb{R}) \rtimes \langle \varphi ^*
\rangle$.

Any isometry $\beta \in \mathbf{\overline{B}_n}(\mathbb{R}) $ of
$\pi ^{-1}(\map)$ is a composition $\tau \ell$ of isometries of
$\pi ^{-1}(\map)$ where $\tau$ acts on $\sigma (\map)$ by
translations and $\ell \in \mathbf{U_n}(\mathbb{R})$. We fix
$\tau$ and $\ell$ in the four paragraphs below and define two
bilipschitz maps $X \rt X$: $H_\infty(\tau)$ and $H_\infty(\ell
)$. Then we will compose these two functions to define a
bilipschitz map $H_\infty(\beta ):X \rt X$.

Right multiplication in a Lie group with a left-invariant metric
is bilipschitz and a finite distance in the sup-norm from the
identity. Therefore, for $u \in \mathbf{U_n}(\mathbb{R})$ and $x$
in a given symmetric apartment $\mathcal{A}$, the map
$(x,u)\mapsto (\tau^{-1}x,\tau u \tau ^{-1})$ defines a
bilipschitz map $\pi ^{-1} (\mathcal{A}) \rt \pi ^{-1}
(\mathcal{A})$ that is a bounded distance from the identity; it is
precisely right multiplication by $\tau ^{-1}$. Thus, we can
compose with the isometry $\tau$ to obtain the bilipschitz map
$\pi ^{-1}(\mathcal{A} ) \rt \pi ^{-1}(\mathcal{A} )$ given by
$(x,n)\mapsto (x,\tau n \tau ^{-1})$. This map is a finite
distance in the sup-norm from $\tau$.

Now we define $H_\infty(\tau ) :X \rt X$ as the map
$$H_\infty(\tau)(x,n)=(x,\tau n \tau ^{-1})$$ From the above paragraph we
know that $H_\infty(\tau)$ restricts to a bilipschitz map on any
symmetric pre-image. It follows from the fact that any path in $X$
is a union of paths in symmetric pre-images that $H_\infty(\tau)$
is bilipschitz.

It is easier to define $H_\infty(\ell ):X\rt X$ as a bilipschitz
map. Just let
$$H_\infty(\ell)(x,n)=(x,a^{-1}_x \ell a_x n)$$
Then $H_\infty(\ell)$ restricts to an isometry of every symmetric
pre-image and, hence, is an isometry of $X$.

We define $H_\infty(\beta )=H_\infty(\tau ) \circ H_\infty(\ell )$
where $\beta = \tau \ell$ as above. Note that $H_\infty(\beta )$
stabilizes each symmetric pre-image, and when restricted to any
symmetric pre-image, it is a finite distance in the sup-norm from
the isometry $\beta$. This is enough to check that we have defined
a homomorphism  $$H_\infty : \mathbf{\overline{B}_n}(\mathbb{R})
\rt \mathcal{QI}(X)$$

\smallskip \noindent \textbf{Coupling real and $p$-adic maps.}
Define the group homomorphism $$H:
\mathbf{\overline{B}_n}(\mathbb{R}) \times
\mathbf{\overline{B}_n}(\mathbb{Q}_p)  \longrightarrow
\mathcal{QI}(X)$$ as $H_\infty \times H_p$.

\bigskip \noindent \textbf{Lemma 2(a).} \emph{The group $G_p$ is a
subgroup of $\mathcal{QI}(X)$.}

\bigskip

\noindent \textbf{Proof.} First we show that the kernel of $H$ is
trivial.

Any distinct pair of symmetric pre-images have infinite Hausdorff
distance. Therefore, if $H( \beta , \alpha)$ is a finite distance
in the sup-norm from the identity, then $H(\beta , \alpha )$
stabilizes every symmetric pre-image. Hence, $\alpha$ stabilizes
every symmetric apartment in $X_p$, and the only element in
$\mathbf{\overline{B}_n}(\mathbb{Q}_p)$ of this sort is the
identity.

Assuming now that $\alpha =1$, we have that $\pi ^{-1}(\map)$ is
stabilized by $H( \beta , \alpha)$, and thus $\beta$ acts as a
bilipschitz map that is equivalent to the identity. It follows
that $\beta =1$ and $H$ is an injective homomorphism.

The lemma follows by noting that if $\varphi^*_p$ and
$\varphi^*_\infty$ are the geometric realizations of $\varphi^*$
on $X_p$ and $\pi ^{-1}(\map)$ respectively, then $\varphi^*_p$
permutes symmetric apartments, and we may assume both that
$\varphi^*_p$ stabilizes $\map$ and that $\varphi^*_p (x)=
\varphi^*_\infty (\sigma(x))$ for all $x \in \map$. Thus, we can
simultaneously apply $\varphi^*_p$ to the $X_p$-factor of $X$ and
$\varphi^* _\infty$ to the symmetric pre-images to obtain an order
$2$ isometry $\varphi ^*_X: X \rt X$. This isometry is not type
preserving on the boundary of $X_p$ (or on the boundary of $\pi
^{-1} (\map)$), so it is an infinite distance in the sup-norm from
any $H( \beta , \alpha )$.

\hfill $\blacksquare$


\pagebreak

 \noindent \Large \textbf{3. All quasi-isometries of
the model space} \normalsize \bigskip

In Section 2 we proved $G_p \leq \mathcal{QI}(X)$. In this section
we will prove $ \mathcal{QI}(X) \leq G_p$. We will immediately
have a proof of Theorem A since we showed in Section 1 that
 $$\mathcal{QI} \big( \sg \big) \cong \mathcal{QI} \big( X
\big)$$

To begin, we let $\phi : X \rt X$ be a $(K,C)$ quasi-isometry. Our
goal is to show that $\phi$ is a finite distance in the sup-norm
from a quasi-isometry of the form $H( \beta , \alpha )(\varphi
^*_X)^k$ where $\beta \in \mathbf{\overline{B}_n}(\mathbb{R})$,
$\alpha \in \mathbf{\overline{B}_n}(\mathbb{Q}_p)$, and $k \in
\{\,0,1\,\}$.

\medskip \noindent \large \textbf{3.1 An isometry of the base building} \normalsize \medskip

Let $\text{Hd}$ denote the Hausdorff distance between subsets of
$X$. The following lemma identifies the distance between points in
the base of $X$ as the distance between their fibers.

\bigskip \noindent \textbf{Lemma 3.1(a).} \emph{If $x,y \in X_p$, then}
$$d(x,y)=\text{Hd}\big(\pi ^{-1}(x)\,,\, \pi ^{-1}(y)\big)$$

\bigskip

\noindent \textbf{Proof.} Let $u \in \mathbf{U_n}(\mathbb{R})$ be
given, and let $u'=a_x u a_x ^{-1} \in \mathbf{U_n}(\mathbb{R})$.
As explained in the previous section, $H_\infty (u' ) :X \rt X$ is
an isometry.

By Lemma 1.1(c), there is a geodesic $\gamma \se X$ from $(x,1)$
to $(y,1)$ with length $d(x,y)$. Therefore, $H_\infty ( u ' )
\gamma$ is a geodesic from $(x,u)$ to $(y, a_y ^{-1} u' a_y)\in
\pi ^{-1}(y)$. Consequently, $$d(x,y) \geq \text{Hd}\big(\pi
^{-1}(x)\,,\,\pi ^{-1}(y)\big)$$

For the opposite inequality, we show that the closest point in
$\pi ^{-1}(y)$ to $(x,1)$ is $(y,1)$; the lemma will follow from
1.1(c). Indeed, recall from 1.1 that $\pi$ is distance
nonincreasing, so for any $u \in \mathbf{U_n}(\mathbb{R})$
 $$d\big((x,1)\,,\,(y,u)\big)\geq d\big((x,1)\,,\,(y,1)\big) = d(x,y)$$

\hfill $\blacksquare$

\bigskip \noindent \textbf{Lemma 3.1(b).} \emph{The quasi-isometry $\phi $ induces an
 isometry $\phi
_\pi : X_p \rt X_p$.}

\bigskip

\noindent \textbf{Proof.} Define $\phi _\pi : X_p \rt X_p$ by
$$\phi _\pi = \pi \circ \phi \circ \sigma (x)$$

Recall that $\pi$ is distance nonincreasing and that $\sigma $ is
an isometric inclusion. Hence, for any pair $x_1 , x_2 \in X_p$ we
have:
$$d\big(\phi _\pi (x_1 ) \,,\, \phi _\pi (x_2 ) \big) \leq Kd(x_1 , x_2) +C$$
Now we work for the other inequality to show that $\phi _\pi $ is
a quasi-isometry.

Recall that $X$ is a fiber bundle, and that each of the spaces
$X$, $\mathbf{U_n}(\mathbb{R})$, and $X_p$ are uniformly locally
finite and uniformly contractible. Furthermore,
$\mathbf{U_n}(\mathbb{R})$ is a manifold, and the system of
apartments in $X_p$ coarsely separates points in $X_p$. These are
all the conditions we need to apply Whyte's theorem \cite{Wh}
which tells us that there are points $y_1,y_2 \in X_p$ such that
$$\text{Hd}\Big( \phi \big(\pi ^{-1} (x_i)\big) \,,\,
 \pi ^{-1} (y_i) \Big) \leq
A$$ for some constant $A=A(K,C)$.

Because $\pi$ is distance nonincreasing,
$$d\big( \phi _\pi (x_i) \,,\,
 y_i \big) \leq
A$$

Using Lemma 3.1(a) and the two preceding inequalities, we have
\begin{align*} d(x_1,x_2) & = \text{Hd}
\Big( \pi ^{-1} (x_1) \,,\, \pi ^{-1} (x_2) \Big) \\
& \leq K\text{Hd} \Big( \phi \big(\pi ^{-1} (x_1) \big) \,,\, \phi
\big( \pi
 ^{-1} (x_2) \big)
\Big) +KC \\
& \leq K\text{Hd} \Big( \pi^{-1} (y_1)  \,,\, \pi^{-1} (y_2)
\Big) +2KA+KC \\
& = Kd \big( y_1  \,,\, y_2
\big) +2KA+KC \\
& \leq Kd \big( \phi _\pi (x_1)  \,,\, \phi _\pi (x_2) \big)
+4KA+KC \end{align*}

We have shown that $\phi _\pi$ is a quasi-isometry. We may further
assume that $\phi _\pi$ is an isometry by Theorem 1.1.3 of
Kleiner-Leeb \cite{K-L}. Indeed, Kleiner-Leeb show that any
quasi-isometry of a higher rank Euclidean building is a bounded
distance from an isometry.

\hfill $\blacksquare$

As the isometry group of $X_p$ is $\mathbf{PGL_n}(\mathbb{Q}_p)
\rtimes \langle \varphi ^* \rangle$, the map $\phi _\pi$ is
identified with an element of $\mathbf{PGL_n}(\mathbb{Q}_p)
\rtimes \langle \varphi ^* \rangle$.

\medskip \noindent \large \textbf{3.2 The base isometry preserves $\msc$} \normalsize \medskip

The following lemma will show that $\phi _\pi$ is identified with
 an element of $\mathbf{\overline{B}_n}(\mathbb{Q}_p) \rtimes
\langle \varphi ^* \rangle$.

\bigskip \noindent \textbf{Lemma 3.2(a).} \emph{The sectors $\phi
_\pi (\msc)$ and $\msc$ contain a common subsector.}

\bigskip

\noindent \textbf{Proof.} Let $\ap $ be a symmetric apartment
containing a subsector of $\phi _\pi (\msc )$. Then $\phi _\pi
^{-1}(\ap )$ is a symmetric apartment as well.

After possibly replacing $\phi$ with an equivalent quasi-isometry,
$\phi$ restricts to a quasi-isometry of symmetric pre-images
$$\pi ^{-1} \big(\phi _\pi ^{-1}(\ap )\big) \rt \pi ^{-1} \big(\ap \big)$$
Applying Theorem 1.1.3 of \cite{K-L}, we may further assume that
this map is an isometry.

If $\msc $ and $\phi _\pi (\msc )$ do not contain a common
subsector, then there are two Weyl chambers in $\pi ^{-1}\big(\phi
_\pi (\msc ) \big) \se \pi ^{-1} \big(\ap \big)$ that have an
infinite Hausdorff distance between them. Name these Weyl chambers
$\mathfrak{C}_1$ and $\mathfrak{C}_2$.

 Since there is only one Weyl chamber in
$\pi ^{-1} (\msc)$ up to Hausdorff equivalence,

$$\text{Hd}\Big(\phi ^{-1} \big( \mathfrak{C}_1 \big) \,,\,  \phi ^{-1} \big(
 \mathfrak{C}_2 \big)\Big) < \infty$$
Thus, $$\text{Hd}\Big(\phi \circ \phi ^{-1} \big( \mathfrak{C}_1
\big) \,,\, \phi \circ \phi ^{-1} \big(
 \mathfrak{C}_2 \big)\Big) < \infty$$
This is a contradiction.

\hfill $\blacksquare$

\medskip \noindent \large \textbf{3.3 An isometry of symmetric pre-images} \normalsize \medskip

Let $\ap$ be a symmetric apartment. By Lemma 3.2(a), the map
$\widehat{\varrho} _{\map , \msc} \circ \phi $ restricts to a
quasi-isometry $\pi ^{-1} (\ap ) \rt \pi ^{-1} (\map )$. Then by
Theorem 1.1.3 of \cite{K-L}, this restriction is a bounded
distance in the sup-norm from an isometry which we name $\phi
_\ap$.

In fact, all symmetric apartments determine the same element of
$\mathbf{PGL_n}(\mathbb{R}) \rtimes \langle \varphi ^* \rangle$:

\bigskip \noindent \textbf{Lemma 3.3(a).} \emph{For any symmetric apartment $\ap \se X_p$,
 we have} $$\phi_\ap =
\phi_{\map}$$

\bigskip

\noindent \textbf{Proof.} The lemma follows from that fact that
$\pi ^{-1} (\map)$  and $\pi ^{-1} (\ap)$ intersect in an
unbounded set, and that isometries of $X_\infty$ are completely
determined by their restriction to any open subset of $X_\infty$.

\hfill $\blacksquare$

\medskip \noindent \large \textbf{3.4 The symmetric pre-image isometry preserves
 $\mathfrak{C}_\mathbf{B}$} \normalsize \medskip

Denote a sector in $\map$ that is opposite to $\msc$ by $\msc
^{\text{op}}$. Note that $\sigma ( \msc ^{\text{op}} ) \se \pi
^{-1} (\map )$ is the unique Weyl chamber in $\pi ^{-1} ( \msc
^{\text{op}})$ up to finite Hausdorff distance, as is
$\widehat{f}^{-1}(\mathfrak{C}_\mathbf{B})$.

Our final lemma identifies $\phi _{\map}$ as an element of
$\mathbf{\overline{B}_n}(\mathbb{R}) \rtimes \langle \varphi ^*
\rangle$.

\bigskip \noindent \textbf{Lemma 3.4(a).} \emph{The isometry
 $\phi _{\map} : \pi ^{-1} (\map ) \rt \pi ^{-1} (\map)$ fixes the
Hausdorff equivalence class of the Weyl chamber
$\widehat{f}^{-1}(\mathfrak{C}_\mathbf{B})$.}

\bigskip

\noindent \textbf{Proof.} The isometry $\phi _{\map}$ fixes $\pi
^{-1} ( \msc )$ up to finite Hausdorff distance by Lemma 3.2(a).
Thus, $\phi _{\map}$ fixes  $\pi ^{-1} ( \msc ^{\text{op}} )$ up
to finite Hausdorff distance as well.

The lemma follows since $\sigma ( \msc ^{\text{op}})$ is the only
Weyl chamber in  $\pi ^{-1} (\msc ^{\text{op}} )$ up to finite
Hausdorff distance.

\hfill $\blacksquare$

\medskip \noindent \large \textbf{3.5 Proof of Theorem A.} \normalsize \medskip

By considering whether or not they preserve types in the boundary,
it is clear that $\phi _{\map} \in
\mathbf{\overline{B}_n}(\mathbb{R})$ if and only if $\phi _\pi \in
\mathbf{\overline{B}_n}(\mathbb{Q}_p)$. Therefore, $\phi$ is a
finite distance in the sup-norm from $H(\phi _{\map} , \phi _\pi
)(\varphi ^*_X)^k$ where $k \equiv 0 $ (mod $2$) if and only if
$\phi _\pi \in \mathbf{\overline{B}_n}(\mathbb{Q}_p)$.

\hfill $\blacksquare$

\pagebreak

 \noindent \Large \textbf{4. Remarks on the large-scale geometry of lamplighter groups} \normalsize
\bigskip

In this section we will describe a geometric model for lamplighter
groups. The model is known as a special case of a Diestel-Leader
graph. We will race through its definition, leaving some small
claims as exercises.

 Using this
model we will explain how one can arrive at Conjectures H and I.

 The connection between Diestel-Leader graphs and lamplighter groups
 is well-known and has played a key role in several results; see
 e.g. \cite{Bt-W}, \cite{Br-W 1}, \cite{Br-W 2}, and \cite{W}.

\bigskip

\smallskip \noindent \textbf{Generalizing solvable
function-field-arithmetic groups.} Throughout this section we are
motivated by virtual lamplighter groups of the form
$\mathbf{B_2}\big(\mathbb{F}_q[t,t^{-1}]\big)$. They are discrete
subgroups of the group
$$\mathbf{B_2}\Big(\mathbb{F}_q((t))\Big) \times
\mathbf{B_2}\Big(\mathbb{F}_q((t))\Big) \; \leq \;
\mathbf{SL_2}\Big(\mathbb{F}_q((t))\Big) \times
\mathbf{SL_2}\Big(\mathbb{F}_q((t))\Big)$$

\smallskip \noindent \textbf{Laurent series with group coefficients.} Let $G$ be a
 finite group. The group of formal Laurent series in a variable
$t$ with coefficients in $G$ is defined as a group of infinite
formal sums:
$$G((t))=\Big\{\,\sum _{i \in \mathbb{Z}} g_i t^i \mid
 g_i \in G\mbox{ for all }i\mbox{ and }g_i=1\mbox{ for }i \ll 0 \,\Big\}$$
Multiplication in $G((t))$ is given by components, so that

$$\Big(\sum  g_i t^i \Big)
\Big( \sum  h_i t^i\Big) = \sum  g_i h_i t^i
$$

\smallskip \noindent \textbf{Laurent polynomials with group
coefficients.} The group of Laurent series contains the subgroup
of Laurent polynomials $$G[t,t^{-1}] =G((t)) \, \cap \,
G((t^{-1}))$$ Notice that the group $\bigoplus _{\mathbb{Z}}G$ is
isomorphic to $G[t,t^{-1}] $.

The search for a space that the lamplighter group $\bigoplus
_{\mathbb{Z}}G \rtimes \mathbb{Z}$ acts on geometrically begins by
finding a space for the group $G((t)) \rtimes \mathbb{Z}$ to act
on.

\smallskip \noindent \textbf{The metric for Laurent series.}
A bi-invariant ultrametric metric on $G((t))$ is given by the norm
$$\Big|\sum g_i t^i\Big| = e ^{-N}$$ where $N$ is
the least integer with $g_N \neq 1$. Thus, the distance between
the two group elements $a,b \in G((t))$ is $|a^{-1}b|$.

\smallskip \noindent \textbf{The regular tree $T_G$.}
We will now construct a regular $(|G|+1)$-valent tree, $T_G$, that
$G((t)) \rtimes \mathbb{Z}$ acts on. The group of Laurent series
will act as ``unipotents" and the integers will act through
hyperbolic isometries.

We pick a line $\ell _a$ for every $a \in G((t))$ along with a
homeomorphism $\rho _a : \mathbb{R} \rt \ell _a$, and we define
the quotient
$$T_G = \Big( \bigcup _{a \in G((t))} \ell _a \Big) \,\Big/\, \sim
$$ where $\rho _a (t) \sim \rho _b (t)$ if $t \geq \log
|a^{-1}b|$. Notice that $T_G$ naturally has the structure of a
simplicial tree with a unique end at infinity corresponding to
each element of the set $G((t)) \cup \{ \infty \}$. We transform
$T_G$ into a metric space by assigning each edge length $1$.

\smallskip \noindent \textbf{Height function on $T_G$.}
There is a well-defined ``height function" $h_G:T_G \rt
\mathbb{R}$ given by $h_G\big( \rho _a (t) \big)=t$ for all $a \in
G((t))$ and $t \in \mathbb{R}$.

\smallskip \noindent \textbf{Action of $G((t)) \rtimes \mathbb{Z}$
 on $T_G$.} By defining $a \cdot \rho _b(t) = \rho _{ab}(t)$, the group
$G((t))$ acts by simplicial automorphisms on $T_G$ that fix the
end corresponding to $\infty$. The height function is invariant
under this action.

 The group $\mathbb{Z}$ acts on $G((t))$ by $n
\cdot \sum g_i t^i = \sum  g_i t^{i-n}$. Therefore, any $n \in
\mathbb{Z}$ acts on $T_G$ by $n \cdot \rho _a (t) = \rho _ {n\cdot
a}(t+n)$. Note that a positive integer $n>0$ expands $G((t))$ and
acts as a hyperbolic isometry that translates its axis, $\ell _1$,
towards the boundary point $\infty$.

Combined, the actions above yield a cocompact but non-proper
action of $G((t)) \rtimes \mathbb{Z}$ on $T_G$ by $(a,n)\cdot
x=a\cdot (n\cdot x )$ for all $x \in T_G$.

\smallskip \noindent \textbf{Action of $G \wr \mathbb{Z}$
 on $T_G \times T_G$.}
There is an automorphism, $\sigma$, of the lamplighter group $G
\wr \mathbb{Z} \cong
 G [t,t^{-1}] \rtimes \mathbb{Z}$ that is defined by
 $$\sigma \Big(\sum  g_i t^i \,,\,m \Big) = \Big( \sum
 g_i t^{-i} \,,\, -m \Big)$$ The automorphism $\sigma$ is used to
include $G \wr \mathbb{Z}$ as a subgroup of $\big(G((t)) \rtimes
\mathbb{Z}\big)
  \times \big(G((t)) \rtimes \mathbb{Z}\big)$ by
$$c \mapsto \big(c\,,\,\sigma (c) \big)$$
This inclusion describes a proper action of $G \wr \mathbb{Z}$ on
$T_G \times T_G$.

\smallskip \noindent \textbf{The horosphere of interest.}
We define the horosphere $$\mathcal{H}_G = \{\,(x,y) \in T_G
\times T_G \mid h_G(x)+h_G(y)=0 \,\}$$ It is easy to check that
$\mathcal{H}_G$ is a locally-finite, connected graph whose path
metric is quasi-isometric to the metric restricted from $T_G
\times T_G$. This graph is a \emph{Diestel-Leader graph}.

The action of $G \wr \mathbb{Z}$ on $T_G \times T_G$ stabilizes
$\mathcal{H}_G$, yielding an action of $G \wr \mathbb{Z}$ on
$\mathcal{H}_G$. The latter action is proper and cocompact.
Therefore, $G \wr \mathbb{Z} $ is quasi-isometric to $
\mathcal{H}_G$.

\smallskip \noindent \textbf{Wreathing infinite groups with $\mathbb{Z}$.}
The definition of $\mathcal{H}_G$ and of the action of $G \wr
\mathbb{Z}$ on $\mathcal{H}_G$ does not require $G$ to be finite.
Our assumption that $G$ is finite is necessary, however, for the
graph $\mathcal{H}_G$ to be locally-finite, and thus for $G \wr
\mathbb{Z} $ to be quasi-isometric to $ \mathcal{H}_G$

\smallskip \noindent \textbf{Aside on finiteness properties.}
It is easy to see that $ \mathcal{H}_G$ contains loops of
unbounded diameter. Thus, the topology of $\mathcal{H}_G$ reflects
that $G \wr \mathbb{Z}$ is finitely generated but not finitely
presented. The former result is trivial; the latter was originally
proved in much greater generality by Baumslag using combinatorial
methods \cite{Ba}.

\bigskip

\smallskip \noindent \textbf{Explaining Conjecture H.}
The table below refers to three basic examples of rank one
amenable groups. The first column denotes the group. Each of these
groups, $\Gamma$, acts properly on a product of negatively-curved
rank one spaces, denoted $X_\Gamma$. The second column lists these
metric spaces. (Below, $T_m$ is an $m+1$-valent tree.)

The third column lists certain subgroups of $\mathcal{QI}(X_
\Gamma )$. Specifically, each group $\Gamma$ is quasi-isometric to
a horosphere $\mathcal{H} \se X_\Gamma$, and
$\mathcal{QI}(X_\Gamma )_\mathcal{H}$ denotes the subgroup of
$\mathcal{QI}(X_\Gamma )$ consisting of classes of
quasi-isometries that preserve $\mathcal{H}$ up to finite
Hausdorff distance.

\bigskip

\begin{center}
\begin{tabular}{llll}  \hline
 \rule[-8pt]{0pt}{22pt} $\Gamma$ \hspace*{.8in} & $X_\Gamma$ \hspace*{.8in} &
  $\mathcal{QI}(X_\Gamma )_\mathcal{H}$
 \\ \hline
\rule{-6pt}{22pt} $\text{BS}(1,m)$ & $\mathbb{H}^2 \times T_m$ &
$\text{Bilip}(\mathbb{R})
\times \text{Bilip}(\mathbb{Q}_m)$ \\
\rule{-6pt}{22pt}  Sol & $\mathbb{H}^2 \times \mathbb{H}^2$ &
$\big[\text{Bilip}(\mathbb{R})
\times \text{Bilip}(\mathbb{R})\big] \rtimes \mathbb{Z} / 2\mathbb{Z}$ \\
\rule{-6pt}{22pt}  $G \wr \mathbb{Z}$ & $T_G \times T_G$ &
$\big[\text{Bilip}\big(G((t))\big) \times
\text{Bilip}\big(G((t))\big)\big] \rtimes \mathbb{Z} /
2\mathbb{Z}$ \vspace*{.1in}
\\  \hline
\end{tabular}
\end{center}

\bigskip

\bigskip

The groups of bilipschitz homeomorphisms that appear in the
right-hand column of the table act on the boundary of a factor of
$X_\Gamma$. Any map on the boundary of one of these rank one
spaces that fixes a distinguished point at infinity, and restricts
to the rest of the boundary points---$\mathbb{R}$ in the case of
$\mathbb{H}^2$ and the $m$-adic numbers in the case of a
$m+1$-valent tree---as a bilipschitz homeomorphism, can easily be
seen to extend to a quasi-isometry of the interior of the space.
(Note that $G((t))$ is isometric to the $|G|$-adic numbers
$\mathbb{Q}_{|G|}$.)

 It is easy to see in all three cases for $\Gamma$, that
$\mathcal{QI}(X_\Gamma )_\mathcal{H} \leq \mathcal{QI}(\Gamma)$.
Theorem B states that the only quasi-isometries of $\Gamma
=\text{BS}(1,m)$ are those accounted for in $\mathcal{QI}(X_\Gamma
)_\mathcal{H}$. That is, that the inclusion above is an equality
for the solvable Baumslag-Solitar groups. Conjectures E and H ask
that the same equality hold for the groups Sol and $G \wr
\mathbb{Z}$ respectively.

\bigskip

\smallskip \noindent \textbf{Explaining Conjecture I.} It is
well-known and easy to see that if $G$ and $H$ are finite groups,
and if $|G|^k=|H|^j$ for some $k,j \in \mathbb{N}$, then $G \wr
\mathbb{Z}$ is quasi-isometric to $H \wr \mathbb{Z}$.
Geometrically, there is a quasi-isometry $T_G \rt T_{G^k}$ defined
by collapsing to a point any path in $T_G$ that bijects under the
height function onto an interval of the form $[ (\ell -1) k \,,\,
\ell k-1 ]$ for $\ell \in \mathbb{Z}$. Applying the above
quasi-isometry to each factor defines a quasi-isometry $T_G \times
T_G \rt T_{G^k} \times T_{G^k}$ that maps $\mathcal{H}_G$ within a
finite Hausdorff distance of $\mathcal{H}_{G^k}$. Similarly,
$\mathcal{H}_H$ and $\mathcal{H}_{H^j}$ are quasi-isometric. It
follows that $|G|^k=|H|^j$ for some $k,j \in \mathbb{N}$ implies
$G \wr \mathbb{Z}$ is quasi-isometric to $H \wr \mathbb{Z}$ since
$\mathcal{H}_{G^k}$ and $\mathcal{H}_{H^j}$ are isometric given
the hypothesis.

There aren't any known examples of a pair of lamplighter groups
that are quasi-isometric but do not meet the hypothesis from the
above paragraph. Given the analogy between lamplighter groups and
solvable Baumslag-Solitar groups, it is possible that such
nonelementary examples do not exist.

Indeed, Farb-Mosher showed that any quasi-isometry between two
solvable Baumslag-Solitar groups---say $\text{BS}(1,m)$ and
$\text{BS}(1,r)$---induce a bilipschitz homeomorphism between the
``upper boundaries" of the respective groups---in concrete terms,
a bilipschitz homeomorphism between $\mathbb{Q}_m$ and
$\mathbb{Q}_r$ (Theorem 6.1 \cite{F-M 1}). Cooper showed that such
a bilipschitz map can only exist when $m^k=r^j$ for some $k,j \in
\mathbb{N}$ (Corollary 10.11 \cite{F-M 1}). If a full analogue
held for lamplighter groups, then a quasi-isometry $G \wr
\mathbb{Z} \rt H \wr \mathbb{Z}$ would induce a bilipschitz
homeomorphism $G((t)) \rt H((t))$ which would imply that
$|G|^k=|H|^j$ for some $k,j \in \mathbb{N}$.

\bigskip

\smallskip \noindent \textbf{A conjecture of Diestel-Leader.} Let
 ${\mathbb{F}_2}$ and ${\mathbb{F}_3}$ be fields with exactly $2$ and $3$
  elements respectively. We define the graph
$$\mathcal{H}_{2,3} = \{\,(x,y) \in T_{\mathbb{F}_2}
\times T_{\mathbb{F}_3} \mid
h_{\mathbb{F}_2}(x)+h_{\mathbb{F}_3}(y)=0 \,\}$$ This graph is
connected, transitive, and locally-finite.

In response to a question of Woess, Diestel-Leader conjectured
that $\mathcal{H}_{2,3}$ is a connected, transitive, and
locally-finite graph that is not quasi-isometric to a
finitely-generated group \cite{D-L}.

Following the reasoning above, the quasi-isometry group of
$\mathcal{H}_{2,3}$ might be $\text{Bilip}(\mathbb{Q}_2) \times
\text{Bilip}(\mathbb{Q}_3)$. Assuming this was the case, if a
finitely-generated group $\Gamma$ was quasi-isometric to
$\mathcal{H}_{2,3}$, then $\Gamma$ would quasi-act on
$\mathcal{H}_{2,3}$, defining a representation $$\Gamma \rt
\text{Bilip}(\mathbb{Q}_2) \times \text{Bilip}(\mathbb{Q}_3)$$
This would make the Diestel-Leader conjecture susceptible to some
algebraic techniques. Perhaps techniques similar to those in
\cite{F-M tukia} could show that no such $\Gamma$ exists.

\pagebreak

\end{document}